\input amstex
\documentstyle{amsppt}
\magnification 1200
\vcorrection{-1cm}
\NoBlackBoxes
\input epsf

\topmatter
\title
       On the hyperbolicity locus of a real curve
\endtitle
\author S.~Yu.~Orevkov
\endauthor
\abstract
Given a real algebraic curve in the projective 3-space,
its hyperbolicity locus is the set
of lines with respect to which the curve is hyperbolic.
We give an example of a smooth irreducible curve whose hyperbolicity locus is
disconnected but the connected components are not distinguished by the
linking numbers with the connected components of the curve.
\endabstract

\address
Steklov Mathematical Institute, Gubkina 8, Moscow, Russia
\endaddress

\address
IMT, l'universit\'e Paul Sabatier, 118 route de Narbonne, Toulouse, France
\endaddress

\email
orevkov\@math.ups-tlse.fr
\endemail

\endtopmatter

\def\refB  {1}
\def\refKS {2}
\def\refMO {3}
\def\refO  {4}
\def\refSV {5}

\def\figXY {1}
\def\figXZ {2}
\def\figPQ {3}

\def\RP{\Bbb{RP}}
\def\LK{\text{lk}} 

\document
\head 1. Introduction \endhead

Let $C$ be an irreducible real algebraic curve in $\RP^n$.
Following [\refKS], we define the {\it hyperbolicity locus} $\Cal H(C)$
as the set of all codimention 2 linear projective subspaces $L$ of $\RP^n$
such that $C$ is {\it hyperbolic} with respect to $L$, i.~e.,
$C$ is disjoint from $L$ and each hyperplane passing throught $L$
has only real intersections with (the complexification of) $C$.

Shamovich and Vinnikov asked [\refSV, Question 3.13] whether $\Cal H(C)$
is always connected. This is evidently true for $n=2$.
However, Kummer and Shaw [\refKS] have shown that the answer is negative for $n=3$.
They gave an example of a sextic genus one smooth real curve $C$ in $\RP^3$
consisting of two topological circles $A$ and $B$ such that
$C$ is hyperbolic with respect to each of some two lines $L$, $L'$
but the linking numbers of $A$ and $B$ with $L$ and $L'$ are different:
$\big(\LK(A,L),\LK(B,L)\big)\ne\big(\LK(A,L'),\LK(B,L')\big)$.

Note that in [\refMO, Theorem 3 and Lemma 3.12] we gave (without stating
it explicitly) an infinite series of examples with any number of
connected components of $\Cal H(C)$. These are the curves
$W_g(a_0,\dots,a_g)$ in the notation of [\refMO], in particular,
the curve constructed in [\refKS] is our $W_1(2,2)$. As in [\refKS], in all
these examples any two  components of $\Cal H(C)$ are distinguished by the
linking numbers. This fact follows from [\refMO, Proposition 3.13].

So, a new question naturally arises ([\refKS, Question 2]): is it possible
that $\Cal H(C)$ is disconnected but the elements of different components
have the same linking numbers with all components of $C$?
Here we give an affirmative answer to this question for $n=3$.

We construct a rational curve $C$ in $\RP^3$ of degree $8$ and two lines
$L,L'\in\Cal H(C)$ which belong to different connected components of $\Cal H(C)$ because
the links $C\cup L$ and $C\cup L'$
are not isotopic in $\RP^3$. The curve $C$ has only one connected component,
thus the linking numbers cannot distinguish between the components of $\Cal H(C)$.

\head 2. The example \endhead

\subhead 2.1. An auxiliary line arrangement
\endsubhead
Let $(x,y,z)$ be coordinates in an affine chart of $\RP^3$.
Let $L$ be the $z$-axis and $L'$ be the common line at infinity of
the planes $z=\text{const}$.
Let $R$ be the rotation by $90^{\circ}$ around $L$. 
We set:
$p_0=(3,-1,-1)$, $q_0=(3,1,1)$, and $p_k=R^k(p_0)$, $q_k=R^k(q_0)$,
$\ell_k=(p_k q_k)$, $\ell'_k=(p_k, q_{k+1})$; see Figure \figXY\
where 
the black and grey parts of
the lines $\ell_k$ and $\ell'_k$ represent the sign of the $z$-coordinate on them.

\midinsert
\epsfxsize=55mm
\centerline{\epsfbox{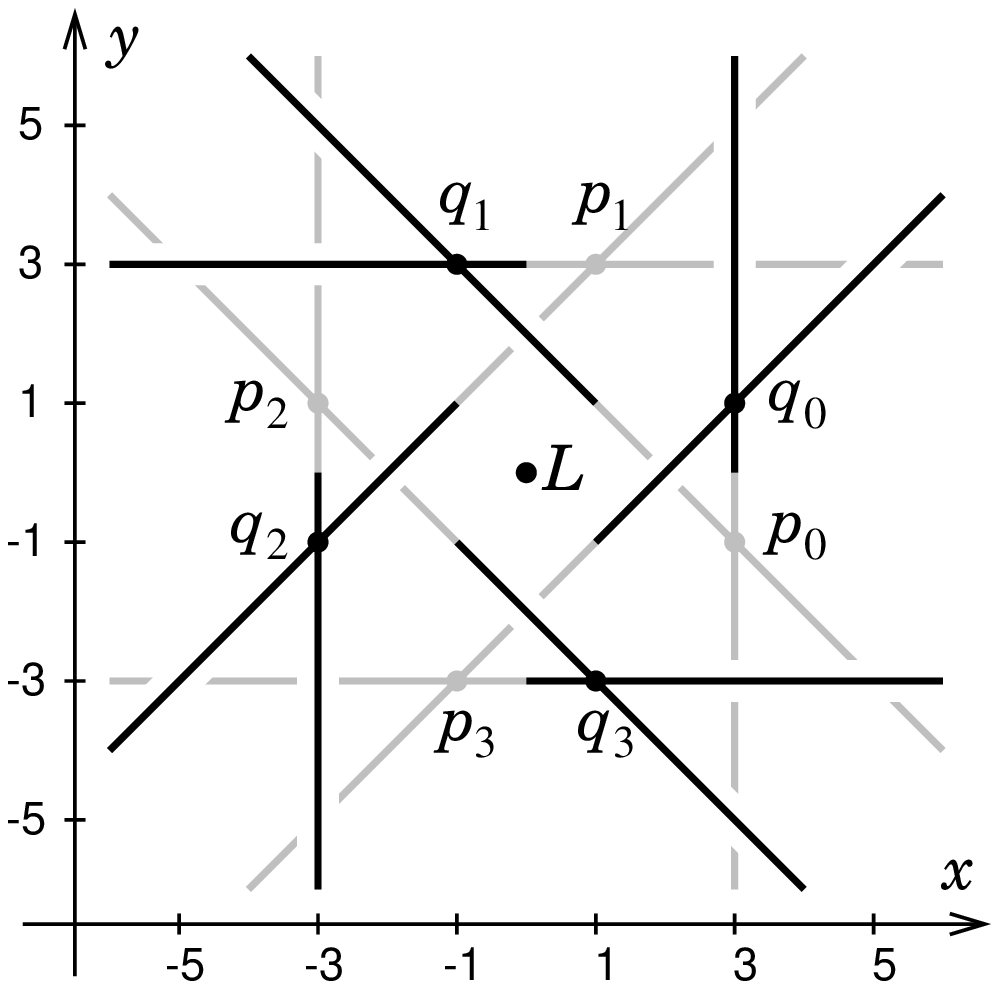}}
\botcaption{Figure \figXY} Projection to $Oxy$ \endcaption
\endinsert

In Figure \figXZ\ we show the same line arrangement in the $Oxz$ projection.
Here the color represents the sign of the $y$-coordinate. Note that since
Figure \figXZ\ is obtained from Figure \figXY\ by a rotation around the axis $Ox$,
the direction of the $y$-axes in Figure \figXZ\ is opposite to
the direction of the $z$-axes in Figure \figXY.

\midinsert
\epsfxsize=55mm
\centerline{\epsfbox{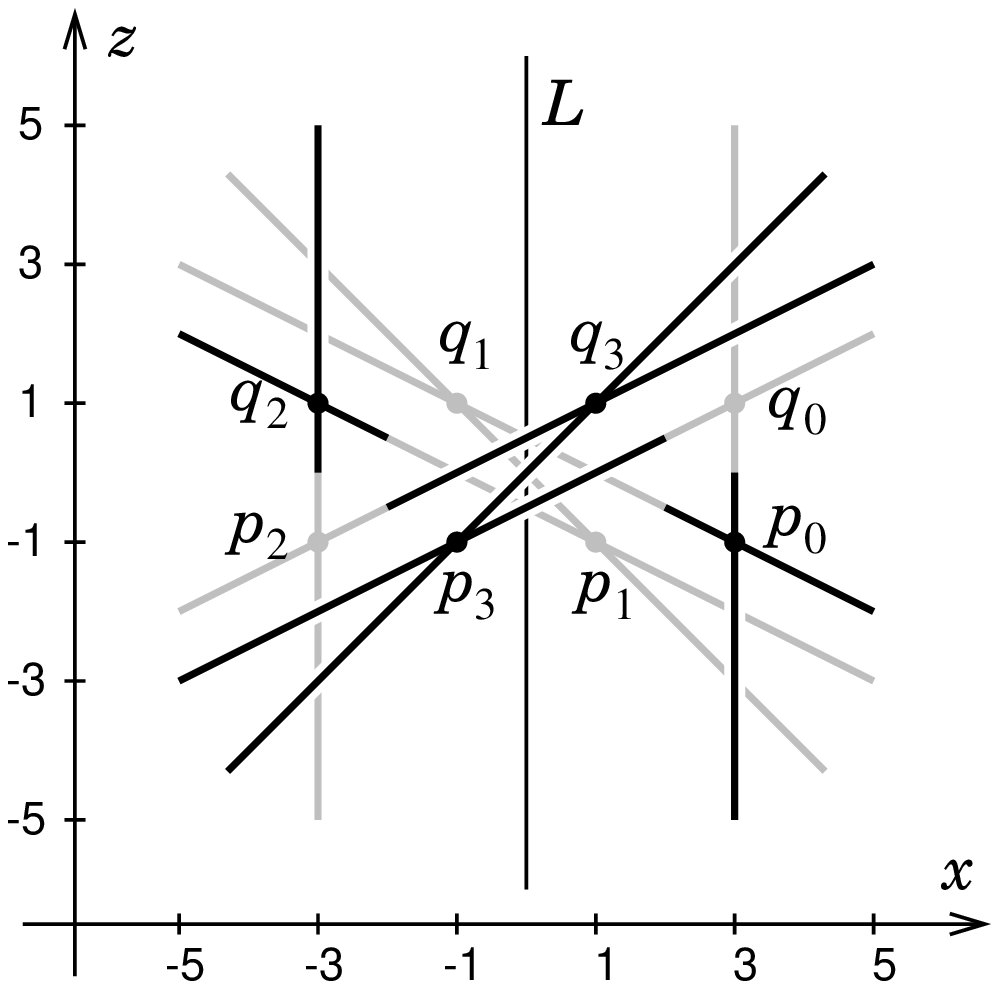}}
\botcaption{Figure \figXZ} Projection to $Oxy$ \endcaption
\endinsert

We orient the lines $\ell_k$ and $\ell'_k$ so that $dz>0$ on them.
Note that in this case we have $d\theta>0$ on them as well
where $(x,y)=(r\cos\theta,r\sin\theta)$.
This means that (under a suitable choice of the orientations of $L$ and $L'$)
both $L$ and $L'$ are positively linked with any of the lines $\ell_k$, $\ell'_k$.

\subhead 2.2. Construction of the curve $C$
\endsubhead
We perturb the union of 8 lines constructed in \S2.1 so that the double point
$q_0$ transforms as in Figure \figPQ(left) and all the other double points ($q_1$, $q_2$,
$q_3$, and $p_0,\dots,p_3$) as in Figure \figPQ(right).
Such a perturbation is possible due to [\refB, Theorem 2.4] (see also [\refMO, Lemma 5.1]):
one should add the lines one by one.
It is easy to see that $L,L'\in\Cal H(C)$.

\midinsert
\epsfxsize=70mm
\centerline{\epsfbox{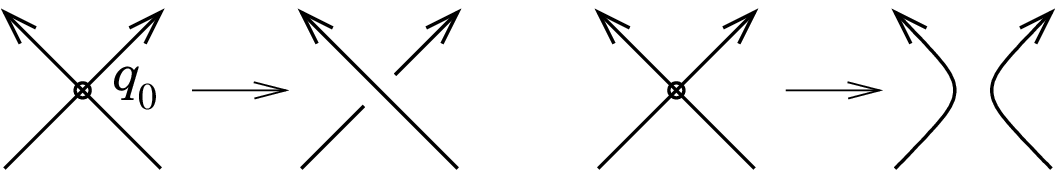}}
\botcaption{Figure \figPQ} Pertubations at $q_0$ (left) and at other points (right)\endcaption
\endinsert

\subhead 2.3. Proof that $L$ and $L'$ belong to different components of $\Cal H(C)$
\endsubhead
It is enough to show that the oriented links $C\cup L$ and
$C\cup L'$ are not isotopic in $\RP^3$. Indeed, their
lifts on the 3-sphere (we denote them $\Lambda$ and $\Lambda'$)
are distinguished by the {\it link determinant}, i.e., the determinant
of the symmetrized Seifert matrix (the value at $-1$
of the Alexander polynomial). Namely, we have $\det(\Lambda)=0$
and $\det(\Lambda')=64$.

For the computations we used the program from [\refO, Appendix]
available at
\chardef\ti=`\~
{\tt http://picard.ups-tlse.fr/\ti orevkov/sm.mat}.
It takes the input in the form of braid. So, for the reader's convenience
we give here the braids $\beta$ and $\beta'$ whose braid closures are
$\Lambda$ and $\Lambda'$ respectively.

To find $\beta'$, we rotate a line by $360^\circ$
around the origin of the plane in
Figure~\figXY\ and write down the contributions of all crossings
consecutively scanned by this line (including the crossings at infinity).
So, we have $\beta' = \beta'_{1/2}\tau_9(\beta'_{1/2})$ where
$\beta'_{1/2}$ is the contribution of the rotation by $180^\circ$ starting
at a horizontal position,
and $\tau_n:B_n\to B_n$ is the braid group automorphism
given by $\sigma_i\mapsto\sigma_{n-i}$.
We need to apply $\tau_9$ on the second half-turn because
the orientation of the line reverses (see also [\refMO, \S\S 4.3--4.5]).
We have
$$
\xalignat 2
\beta'_{1/2}=\,
  &\sigma_1 \Delta_{45} \sigma_8      && \pm(2,0),\;(\infty,0)\\
  &\times\sigma_2^{-1}                && (3,1)\;\;\text{(no contrib. of $(-3,-1)$)}\\
  &\times\sigma_3 \sigma_6            && \pm(5,3)\\
  &\times\sigma_2\Delta_{45}\sigma_7  && \pm(3,3),\;(\infty,\infty)\\
  &\times\sigma_3 \sigma_6.           && \pm(3,5)\\
  &\sigma_1 \Delta_{45} \sigma_8      && \pm(0,2),\;(0,\infty)\\
  &\times\sigma_3 \sigma_6            && \pm(-3,5)\\
  &\times\sigma_2\Delta_{45}\sigma_7  && \pm(-3,3),\;(-\infty,\infty)\\
  &\times\sigma_3 \sigma_6            && \pm(-5,3)
\endxalignat
$$
where $\Delta_{45}=\sigma_4\sigma_5\sigma_4$ is the contribution of each
triple crossing at infinity (in the comments we refer to the coordinates
of the contributing crossings in Figure \figXY).

Similarly, $\beta=\beta_{1/2}\tau_9(\beta_{1/2})$ with
$$
\split
   \beta_{1/2} =& 
             (\sigma_1\sigma_4\sigma_7)
             \sigma_2^{-1}
             (\sigma_3 \sigma_5)
             (\sigma_2\sigma_4\sigma_6)
             (\sigma_3 \sigma_5)
           \\
            &\times
             (\sigma_1\sigma_4\sigma_7)
             (\sigma_3 \sigma_5)
             (\sigma_2\sigma_4\sigma_6)
             (\sigma_3 \sigma_5)
           \\
           &\times
           (\sigma_8\sigma_7\dots\sigma_1)
           \endsplit
$$
where the first two lines represent $C$
and the third line represents $L'$.

\medskip
{\bf Remark.} If we replace the negative crossing
in Figure \figPQ(left) by a positive one
(which corresponds to replacing $\sigma_2^{-1}$ by $\sigma_2$ in our braids),
then $\Cal H(C)$ will be connected by [\refMO, Proposition 3.13] because the obtained
8-th degree curve will be maximally writhed in this case.


\enddocument